\documentclass[10pt]{article}
\usepackage{amsfonts}
\usepackage{amsmath,amsthm,amssymb,latexsym,amscd}
\usepackage{epsf}
\usepackage[pdftex]{graphicx}
\usepackage{tikz}
\theoremstyle{definition}

\theoremstyle{remark}

\numberwithin{equation}{section}

\usepackage{tikz}
\newcommand{\cG}{{\mathcal{G}}}

\newcommand{\cH}{{\mathcal{H}}}

\newcommand{\cW}{{\mathcal{W}}}

\newcommand{\cJ}{{\cal J}}

\newcommand{\cS}{{\mathcal S}}

\newcommand{\cA}{{\cal A}}

\newcommand{\eps}{{\epsilon}}

\title{\bf Complex fracture nucleation and evolution with nonlocal elastodynamics }

\author{Robert Lipton\thanks{Louisiana State University, Baton Rouge, USA}
\and
Stewart Silling
\thanks{Sandia National Laboratories, Albuquerque, USA}
\and
Richard Lehoucq
\thanks{Sandia National Laboratories, Albuquerque, USA}}

\date{}
\begin{document}
\maketitle
\begin{abstract}
{A mechanical model is introduced for predicting the initiation and evolution of complex fracture patterns without the need for a damage variable or law. 
The model, a continuum variant of Newton's second law, uses integral rather than partial differential operators where the region of integration is over finite domain. 
The force interaction is derived from a novel nonconvex strain energy density function, resulting in a nonmonotonic material model.
The resulting equation of motion is proved to be mathematically well-posed.
The model has the capacity to simulate nucleation and growth of multiple, mutually interacting dynamic fractures.
In the limit of zero region of integration, the model reproduces the classic Griffith model of brittle fracture.
The simplicity of the formulation avoids the need for supplemental 
kinetic relations that dictate crack growth or the need for an 
explicit damage evolution law.}
\end{abstract}
\begin{flushleft}
\end{flushleft}
%\abbreviations{SAM, self-assembled monolayer; OTS, octadecyltrichlorosilane}
\pagestyle{myheadings}
\markboth{}{Complex Fracture}

\section{Introduction}\label{sec-intro}

%\dropcap
{S}imulation of dynamic fracture is a challenging problem 
because of the extremes of strain and strain-rate experienced by the
material near a crack tip, and because of the inherent instabilities such 
as branching that characterize many applications.
These considerations, as well as the incompatibility of  partial 
differential equations (PDEs) with discontinuities,
have led to the formulation of specialized methods for the simulation of 
crack growth, especially in finite element analysis.
These techniques include the extended finite element \cite{moes2002}, \cite{Duarte},
cohesive element \cite{elices2002}, and phase field \cite{bourdin}, \cite{Borden}, \cite{Miehe} methods and have met with notable successes.

%The peridynamic theory of solid mechanics \cite{silling-pdreview} has 
%been proposed to be compatible along discontinuities.
%The peridynamic theory of solid mechanics \cite{silling-pdreview} has been proposed as a generalization of the standard theory of solid
%mechanics by overcoming the aforementioned incompatibility of a PDE model. (RBL Jan 4 2016).
%-that is more compatible than the standard PDEs with the physical nature of cracks as discontinuities.
The peridynamic theory of solid mechanics \cite{silling-pdreview} has been proposed as a generalization of the standard theory of solid mechanics that predicts the creation and growth of cracks. In this formulation crack dynamics is given directly by evolution equations for the deformation field eliminating the need for supplemental kinetic relations describing crack growth.  %(RPLJan 11 2016)
The balance of linear momentum takes the form
\begin{equation}
   \rho(x)u_{tt}(x,t) = \int_{\cH_\eps(x)} f(y,x)\;dy + b(x,t)
\label{eqn-uf}
\end{equation}
where $\cH_\eps(x)$ is a neighborhood of $x$,  $\rho$ 
is the density, $u$ is the displacement field, $b$ is the body force 
density field, and $f$ is a material-dependent
function that represents the force density (per unit volume squared) 
that point $y$ exerts on $x$ as a result of the deformation.
The radius $\eps$ of the neighborhood is referred to as the \emph{horizon}. 
%This formulation avoids the mathematical difficulty of trying to apply a PDE (in either the strong or the weak form)  on a singularity such as a crack.
The motivation for peridynamics is that all material points are subject to the same basic field equations, whether on or off of a discontinuity; the equations also have a basis in non-equilibrium
statistical mechanics  \cite{statmech}.
This paradigm, to the extent that it is successful, liberates analysts from the need to develop and implement supplementary equations that dictate the evolution of discontinuities.

%These nonlocal balance laws~\cite{statmech} can be derived using the principles of non-equilibrium statistical mechanics as expectation in phase space.

Standard practice in peridynamics dictates that the nucleation 
and propagation of cracks requires the specification of a damage variable 
within the functional form of $f$ that irreversibly degrades or eliminates the pairwise force interaction
between $x$ and its neighbor $y$. This is referred to as breaking the {\emph{bond}} between $x$ and $y$. 
Here the term ``bond'' is used only to indicate a force interaction between two material 
points $x$ and $y$ through some potential, whose value can depend on the deformations of
other bonds as well.
A wide variety of damage laws in peridynamics are possible, 
and often they contain parameters that can be calibrated to
important experimental measurements such as critical energy release 
rate \cite{silling05} or the Eshelby-Rice $J$-integral \cite{hu2012}.
Damage evolution in peridynamic mechanics can be cast in a consistent 
thermodynamic framework \cite{silling-pdreview}, including appropriate 
restrictions derived from the Second Law of thermodynamics.
This general approach of using bond damage has met with notable 
successes in the simulation of dynamic fracture \cite{ha2010,foster2011}.
However, because of the large number of bonds in a discrete formulation of \eqref{eqn-uf}, there is a cost associated with keeping track of bond damage,
as well as the need to specify a bond damage evolution law.

In the present paper, we report on recent efforts to model cracks in peridynamics without a bond damage variable. 
The main innovation in the present paper is a nonconvex elastic material model for peridynamic mechanics that, 
under certain conditions, nucleates and evolves discontinuities spontaneously.
This approach is rigorously shown to reproduce the most salient experimentally observed characteristic of brittle fracture---the nearly constant
amount of energy consumed by a crack per unit area of crack growth (the Griffith crack model).
Our results further show that in spite of the strong nonlinearity of the material model, the resulting equation of motion 
is well-posed within a suitable function space, providing a 
mathematical context for which multiple interacting cracks can grow without recourse to supplemental kinetic relations.
In the limit of small horizon $\eps$, the nonconvex peridynamic model recovers a limiting fracture evolution characterized by the classical PDE of linear elasticity away from the cracks. 
The evolving fracture system for the limit dynamics is shown to have bounded Griffith fracture energy described by a critical energy release rate obtained directly from the nonconvex peridynamic  potential. 
These results bring the field of peridynamic mechanics closer to the goal of generalizing the conventional theory to model both continuous and discontinuous 
deformation using the same balance laws. 

\section{Nonconvex material model}

Let $\cS$ denote the {\emph{bond strain}}, defined to be the change in the length of a bond as a result of 
deformation divided by its initial length.
We assume that the displacements $u$ are small (infinitesimal) relative to the size of the body $D$. 
Under this hypothesis the strain between two points $x$ and $y$ under the displacement field $u$ is given by  
\begin{equation}
\cS_u=\frac{u(y,t)-u(x,t)}{|y-x|} \cdot e\,, \qquad e=\frac{y-x}{|y-x|}\,,
\label{strainpd}
\end{equation}
where $e$ is the unit vector in the direction of the bond and $\cdot$ is the dot product between two vectors.
To describe the material response, assume that the force interaction between
points $x$ and $y$ reversibly stores potential (elastic) energy, and that this energy depends
only on the bond strain and the bond's undeformed length.
The elastic energy density at a material point $x$ is assumed to be given by
\begin{equation}
  W(x)=\frac{1}{V_\eps}\int_{\cH_\eps(x)} |y-x|\cW^\eps\big(\cS_u,y-x\big)\;dy
\label{eqn-W}
\end{equation}
where $\cW^\eps(\cS,y-x)$ is the pairwise force potential per unit length between $x$ and $y$ and
$V_\eps$ is the area (in 2D) or the volume (in 3D) of the neighborhood $\cH_\eps(x)$.

The nonconvexity of the potential $\cW^\eps$ with respect to the strain $\cS$ distinguishes this material model from those previously considered in the peridynamic literature.
By Hamilton's principle applied to a bounded body $D\subset\mathbb{R}^d$, $d=2$, $3$,  
the equation of motion describing the displacement field $u(x,t)$ is
\begin{equation}
\rho\, u_{tt}(x,t)=\frac{2}{V_\eps}\int_{\mathcal{H}_\epsilon(x)}\,
  \big(\partial_{\cS} \mathcal{W}^\epsilon(\mathcal{S}_u,y-x)\big)e\,dy+b(x,t),
\label{eqofmotion}
\end{equation}
which is a special case of \eqref{eqn-uf}.
The evolution described by \eqref{eqofmotion} is investigated in detail in the papers \cite{Lipton1,Lipton2}. 

We assume the general form
\begin{equation}
   \cW^\epsilon(\cS,y-x)=\frac{J^\epsilon(|y-x|)}{\epsilon|y-x|} \Psi(|y-x|\cS^2)
\label{eqn-WJPsi}
\end{equation}
where $J^\epsilon(|y-x|)=J(|y-x|/\eps)>0$ is a weight function and
$\Psi:[0,\infty)\rightarrow\mathbb{R}^+$  is a continuously 
differentiable function such that $\Psi(0)=0$, $\Psi'(0)>0$, and $ \Psi_\infty:=\lim_{r\rightarrow
\infty}\Psi(r) <\infty$. 
The pairwise force density is then given by
\begin{equation}
\partial_\cS\mathcal{W}^\epsilon(\cS,y-x)=\frac{2J^\epsilon(|y-x|)}{\epsilon}\Psi'\big(|y-x|\cS^2\big)\cS.
\label{force}
\end{equation}
For fixed $x$ and $y$, there is a unique maximum in the curve of force versus  strain (Figure~1).
The location of this maximum can depend on the distance between $x$ and $y$ and
occurs at the bond strain $\cS_c$ such that $\partial^2\cW^\eps/\partial\cS^2(\cS_c,y-x)=0$.
This value is $\cS_c=\sqrt{r_c/|y-x|}$, where $r_c$ is the unique number such that
$\Psi'(r_c)+2r_c\Psi''(r_c)=0$.

%\begin{figure} % [f]
%\includegraphics[width=.45\textwidth]{fig-nonconvex}
%\caption{{\bf Relation between force and  strain for $x$ and $y$ fixed.}}
%\label{fig-nonconvex}
%\end{figure}

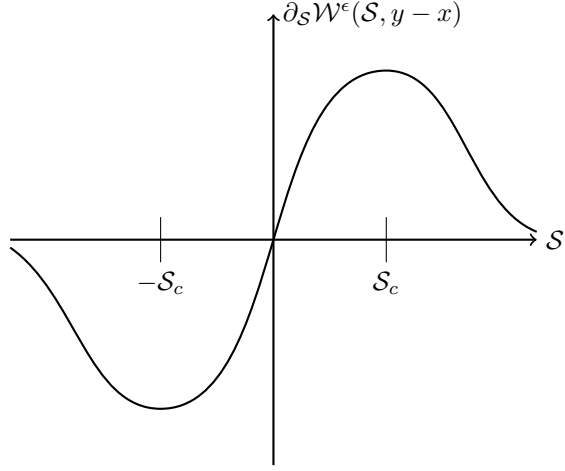
\begin{figure} 
\label{fig-nonconvex}
\centering
\begin{tikzpicture}[xscale=1.0,yscale=1.5]
\draw [<-,thick] (0,2) -- (0,-2);
\draw [->,thick] (-3.5,0) -- (3.5,0);
\draw [-,thick] (-3.5,-0.07) to [out=-25,in=180] (-1.5,-1.5) to [out=0,in=180] (1.5,1.5)
           to [out=0,in=165] (3.5,0.07);
\draw (1.5,-0.2) -- (1.5, 0.2);
\draw (-1.5,-0.2) -- (-1.5, 0.2);
\node [below] at (1.5,-0.2) {$\cS_c$};
\node [below] at (-1.5,-0.2) {$-\cS_c$};
\node [right] at (3.5,0) {$\cS$};
\node [right] at (0,2.0) {$\partial_\cS\cW^\epsilon(\cS,y-x)$};
\end{tikzpicture} 
\caption{{\bf Relation between force and  strain for $x$ and $y$ fixed.}}
\end{figure}
We introduce $Z(x)$, the maximum value of bond strain relative
to the critical strain $\cS_c$ among all bonds connected to $x$:
\begin{equation}
  Z(x)=\max_{y\in\cH_\eps(x)}\frac{\cS_u(x,y)}{\cS_c(x,y)}.
\label{eqn-Sdef}
\end{equation}
The fracture energy $\mathcal{G}$ associated with a crack is stored in the bonds  corresponding to points $x$ for which $Z(x)\gg1$. It is associated with bonds so far out on the the curve in Figure~1 that they sustain negligible force density. This set contains the {\em{jump set}} $\mathcal{J}_{u}$, 
along which the displacement $u$ has jump discontinuities.

Consider an initial value problem for the body $D$
with bounded initial displacement field $u_0$, bounded initial velocity 
field $v_0$, and a non-local Dirichlet condition $u=0$ for $x$ within a layer of thickness
$\epsilon$ external to $D$ containing the domain boundary $\partial D$.  
The initial displacement $u_0$ can contain a jump set $\mathcal{J}_{u_0}$ associated with an initial network of cracks.
%along which $u_0$ has jump discontinuities associated with a crack. 

This initial value problem for \eqref{eqofmotion} is well posed
provided we frame the problem in the space of square integrable displacements satisfying the nonlocal Dirichlet boundary conditions. 
This space is written $L^2_0(D;\mathbb{R}^d)$.
The body force $b(x,t)$ is prescribed for $0\leq t \leq T$ 
and belongs to $C^1([0,T];L^2_0(D;\mathbb{R}^d))$.  
The papers \cite{Lipton1,Lipton2} establish that
if the initial data $u_0$, $v_0$ are in  $L^2_0(D;\mathbb{R}^d)$, and if $u_0$ has bounded total strain energy,
then there exists a unique solution $u(x,t)$ of \eqref{eqofmotion} 
belonging to $C^2([0,T];L^2_0(D;\mathbb{R}^d))$ taking on the intial data $u_0$, $v_0$.  

%For all elastic peridynamic materials, the expected global balance of energy holds \cite{Silling2,Lipton1,Lipton2}:
%\begin{equation}
%  \cE(t)=\cE(0)+\int_0^t\int_{D} b\cdot u_t(\tau)\,dx\,d\tau
%\label{BalanceEnergy}
%\end{equation}
%where
%\begin{equation}
%  \cE(t)=\frac{\rho}{2}\Vert u_t(t)\Vert^2+P(t), \qquad P(t)=\int_D W\;dx
%\label{energyt}
%\end{equation}
%where  $\Vert q\Vert=(\int_D|q|^2 dx)^{1/2}$ and $W$ is given by \eqref{eqn-W}.

\section{Crack nucleation as a material instability}

Normally, we expect an elastic spring to ``harden,'' that is,  force increases with strain.
If instead the spring ``softens'' and the force decreases, then it is unstable: under constant load, its extension
will tend to grow without bound over time.
A material model of the type shown in Figure~1 has this type of softening behavior for sufficiently large strains.
Yet the instability of a bond between a \emph{single} pair of points $x$ and $y$ does not necessarily imply that the entire
body is dynamically unstable.
Here, we present a condition on the material stability with regard to the growth of infinitesimal jumps in displacement across surfaces.

Let $\gamma(x,t)$ denote the volume fraction of points $y\in\cH_\eps(x)$ such that $\cS_u >\cS_c$.\footnote{This can be thought of as the ``number of bonds'' strained past the threshold divided by the total ``number of bonds'' connected to $x$.}
We apply a linear perturbation analysis of \eqref{eqofmotion} to show 
that small scale jump discontinuities in the displacement can become unstable and grow under certain conditions.

Consider a time independent body force density $b$ and a smooth 
solution $u^\ast$ of \eqref{eqofmotion}. 
Let $x$ be a fixed point in $D$.
We investigate the evolution of a small jump in displacement of the 
form
\[
  u(y,t)=u^\ast(y,t)+\left\{\begin{array}{ll}
     0          & {\textrm{if}}\;(y-x)\cdot n<0, \\
     \bar us(t) & {\textrm{otherwise}}.
   \end{array} \right.
\]
where $\bar u$ is a vector, $s(t)$ is a scalar function of time, and $n$ is a unit vector.
Geometrically, the surface of discontinuity passes through $x$ and has normal $n$.
The vector $\bar u$ gives the direction of motion of points on either side of the surface as they separate.

We give conditions for which the jump perturbation is exponentially unstable.
The {\emph{stability tensor}} $\cA_n(x)$ is defined by
\begin{equation}
  \cA_n(x)=\int_{\cH^-_\epsilon(x)}\frac{1}{|y-x|}\partial^2_{\mathcal{S}}\mathcal{W}^\epsilon(\mathcal{S}_{u^\ast},y-x) \, e\otimes e \; dy\,,
\label{instabilitymatrix}
\end{equation}
where $ {\mathcal{H}}_\epsilon^{-}(x) = \{ y\in{\mathcal{H}}_\epsilon(x) \vert (y-x)\cdot n <0 \} $.
%where $\cH^-_\epsilon(x)=\{x'\in\cH_\epsilon(x)\; | \; (x'-x)\cdot n<0\}$. and
%\[  \cS^\ast = \frac{u^\ast(y)-u^\ast(x)}{|y-x|}.   \]
A sufficient condition for the rapid growth of small jump discontinuity is derived in \cite{Lipton1,Lipton2,SillingWecknerBobaru}.
If the stability matrix $\cA_n(x)$ has at least one negative eigenvalue then (1) $\gamma(x)>0$, and (2) there exist a non-null vector $\bar u$ and a unit vector $n$ such that $s(t)$ grows exponentially in time.
%\begin{itemize}
%\item $\gamma(x)>0$, and
%\item there exist a non-null vector $\bar u$ and a unit vector $n$ such that $s(t)$ grows exponentially in time.
%\end{itemize}
The significance of this result is that the nonconvex bond strain energy model can spontaneously nucleate
cracks without the assistance of supplemental criteria for crack nucleation.
This is an advantage over conventional approaches because crack initiation is predicted by
the fundamental equations that govern the motion of material particles.
A negative eigenvalue of $\cA_n(x)$ can occur only if a sufficient fraction of the bonds connected to $x$
have strains $\cS_{u^\ast}>\cS_c$.

\section{Small horizon limit: dynamic fracture}
For finite horizon $\epsilon>0$ the  elastic moduli and critical energy release rate are recovered 
directly from the strain potential $\mathcal{W}^\epsilon(\mathcal{S},y-x)$ given by \eqref{eqn-WJPsi}.
First suppose the displacement inside $\mathcal{H}_\epsilon(x)$ is affine, that is, $u(x)=Fx$ where $F$ is a constant matrix.  
For small strains, i.e., $\mathcal{S}=Fe\cdot e\ll\mathcal{S}_c$, the strain potential is linear elastic to leading order and characterized by elastic moduli $\mu$ and $\lambda$ associated with a linear elastic isotropic material  
\begin{eqnarray}
W(x)&=&\frac{1}{V_d}\int_{H_\epsilon(x)}|y-x|\mathcal{W}^\epsilon(\cS_u,y-x)\,dy\nonumber\\
&=&2\mu |F|^2+\lambda |Tr\{F\}|^2+O(\epsilon|F|^4).
\label{LEFMequality}
\end{eqnarray}
The elastic moduli $\lambda$ and $\mu$ are calculated directly from the strain energy density \eqref{eqn-WJPsi} and are given by
\begin{equation}\label{calibrate1}
\mu=\lambda=M\frac{1}{d+2} \Psi'(0) \,, 
\end{equation}
where the constant $M=\int_{0}^1r^{d}J(r)dr$ for dimensions $d=2,3$.
In regions of discontinuity the same strain potential \eqref{eqn-WJPsi}  is used to calculate the 
amount of energy consumed by a crack per unit area of crack growth, i.e., the critical energy release rate $\cG$.
Calculation applied to \eqref{eqn-WJPsi} shows that $\mathcal{G}$
equals the work necessary to eliminate force interaction on either side of a fracture surface per unit fracture area and is given in
three dimensions by
\begin{equation}
  \mathcal{G}=\frac{4\pi}{V_d}\int_0^\epsilon\int_z^\epsilon\int_0^{\cos^{-1}(z/\zeta)}
   \mathcal{W}^\epsilon(\infty,\zeta)\zeta^2\sin{\phi}\,d\phi\,d\zeta\,dz
\label{calibrate2formula}
\end{equation}
where $\zeta=|y-x|$.
(See Figure~2 for an explanation of this computation.)  
In $d$ dimensions, the result is
\begin{equation}
  \mathcal{G}=M\frac{2\omega_{d-1}}{\omega_d}\, \Psi_\infty\,,
\label{calibrate2formulad}
\end{equation}
where $\omega_{d}$ is the volume of the $d$ dimensional unit ball, $\omega_1=2,\omega_2=\pi,\omega_3=4\pi/3$.
%Identical results also also apply to the two dimensional problem $(d=2)$. For these models the fracture toughness $\mathcal{G}$ is independent of $\epsilon$. 

\begin{figure} % [f]
\centering
\includegraphics[width=.45\textwidth]{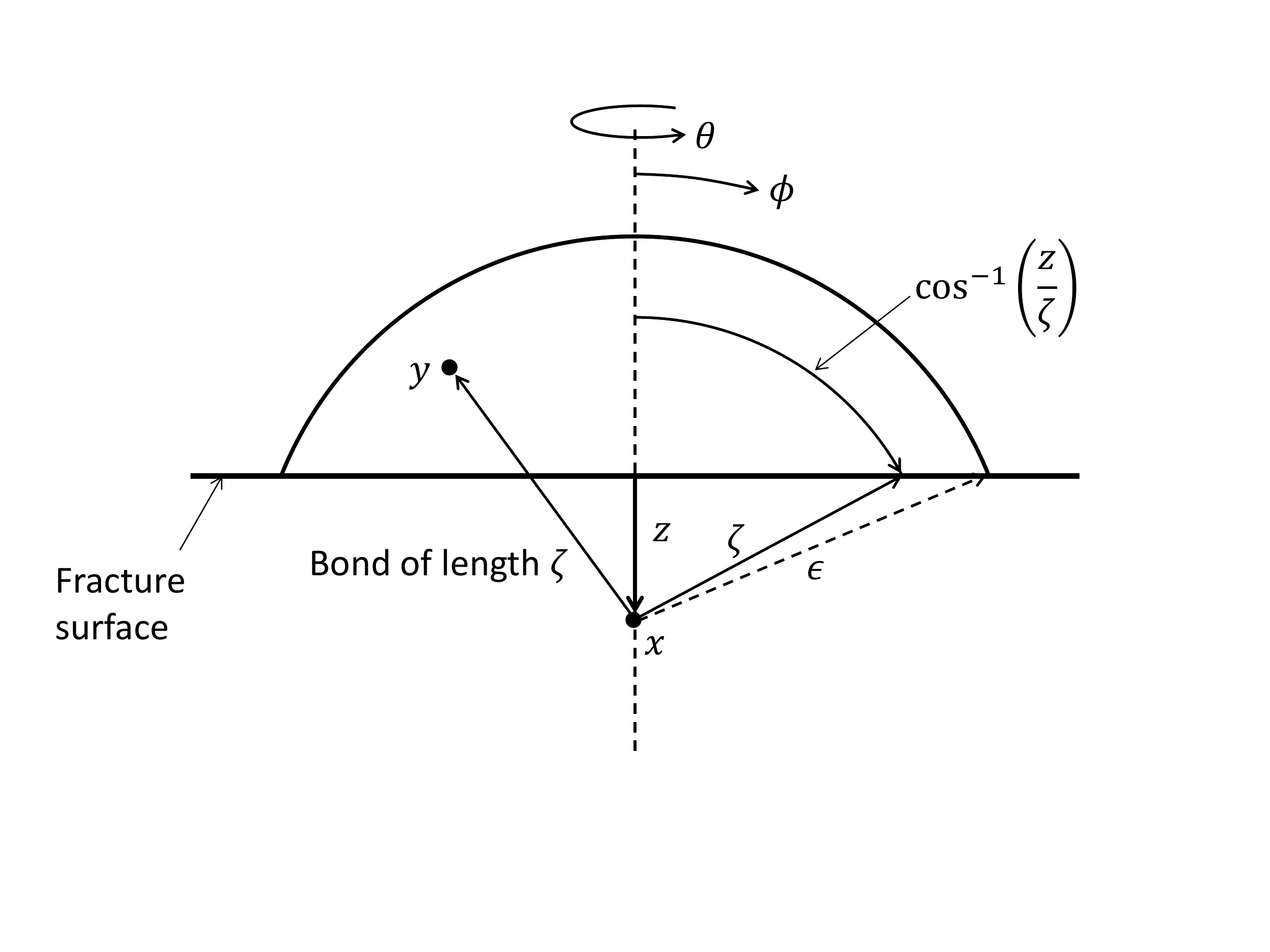}
\label{fig-Gintegral}
\caption{{\bf Evaluation of the critical energy release rate $\mathcal{G}$. 
For each point $x$ along the dashed line, $0\leq z\leq \epsilon$, the work required to break the 
interaction between $x$ and $y$ in the spherical cap is summed up in \eqref{calibrate2formula} using spherical coordinates centered at $x$,
which depends on $z$.}}
\end{figure}

In the limit of small horizon $\epsilon\rightarrow 0$ peridynamic solutions converge in mean square to limit solutions that are 
linear elastodynamic off the crack set, that is, the PDEs of the local theory hold at points off of the crack.
The elastodynamic balance laws are characterized by elastic moduli $\mu$, $\lambda$.
The evolving crack set possesses bounded Griffith surface free energy associated with the critical energy release rate $\mathcal{G}$. 
We prescribe a small  initial displacement field $u_0(x)$ and small initial velocity field $v_0(x)$  with bounded Griffith fracture energy given by
\begin{eqnarray}
\int_{D}\,2\mu |\mathcal{E} u_0|^2+\lambda |{\rm div}\, u_0|^2\,dx+\mathcal{G}|\mathcal{J}_{u_0}|\leq C
\label{LEFMboundinitil}
\end{eqnarray}
for some $C<\infty$.
Here $\mathcal{J}_{u_0}$ is the initial crack set across which the displacement $u_0$ has a jump discontinuity.
This jump set need not be geometrically simple; it can be a complex network of cracks.
$|\mathcal{J}_{u_0}|=H^{d-1}(J_{u_0})$ is the $d-1$ dimensional Hausdorff measure of the jump set.
This agrees with the total surface area (length) of the crack network for sufficently regular cracks for $d=3(2)$.  
The strain tensor associated with the initial displacement $u_0$ is denoted by $\mathcal{E} u_0$. 
Consider the sequence of solutions $u^{\epsilon}$ of the initial value problem associated with progressively smaller peridynamic horizons $\epsilon$.
The  peridynamic evolutions $u^{\epsilon}$   converge in mean square uniformly in time to a 
limit evolution $u^0(x,t)$ in $C([0,T];L_0^2(D,\mathbb{R}^d )$ and $u_t^0(x,t)$ in $L^2([0,T]\times D;\mathbb{R}^{d})$ with the same initial data, i.e.,
\begin{eqnarray}
\lim_{\epsilon\rightarrow 0}\max_{0\leq t\leq T}\int_D|u^{\epsilon}(x,t)-u^0(x,t)|^2\,dx=0,
\label{unifconvg}
\end{eqnarray}
see \cite{Lipton1}, \cite{Lipton2}.
It is found  that the limit evolution $u^0(t,x)$ has bounded Griffith surface energy and elastic energy given by
\begin{eqnarray}
\int_{D}\,2\mu |\mathcal{E} u^0(t)|^2+\lambda |{\rm div}\,u^0(t)|^2\,dx+\mathcal{G}|\mathcal{J}_{u^0(t)}|\leq C, 
\label{LEFMbound}
\end{eqnarray}
for $0\leq t\leq T$,
where $\mathcal{J}_{u^0(t)}$ denotes the evolving fracture  surface inside the domain $D$,  \cite{Lipton1}, \cite{Lipton2}. 
The limit evolution $u^0(t)$ is found to lie in the space of functions of bounded deformation SBD, see \cite{Lipton2}. 
For functions in SBD the bond strain $\mathcal{S}_u$ defined by \eqref{strainpd} is related to the strain tensor $\mathcal{E}u$ by
\begin{equation}
\lim_{\epsilon\rightarrow 0}\frac{1}{V_\epsilon}\int_{\mathcal{H}_\epsilon(x)}|\mathcal{S}_u-\mathcal{E}u(x) e\cdot e|\,dy,
\label{equatesandE}
\end{equation}
for almost every $x$ in $D$.
The jump set $\cJ_{u^0(t)}$  is the countable union of rectifiable surfaces  
(arcs) for $d=3(2)$, see \cite{AmbrosioCosicaDalmaso}.

In domains away from the crack set the limit evolution satisfies local linear elastodynamics (the PDEs of the standard
theory of solid mechanics).
Fix a tolerance $\tau>0$. 
If for subdomains $D'\subset D$ and for times $0<t<T$ the associated strains $\mathcal{S}_{u^\epsilon}$ 
satisfy $|\mathcal{S}_{u^\epsilon}|<S_c$ for every $\epsilon<\tau$  
then it is found that the limit evolution $u^0(t,x)$  is governed by the PDE
\begin{eqnarray}
\rho u^0_{tt}(t,x)= {\rm div}\sigma(t,x)+b(t,x), \hbox{on $[0,T]\times D'$},
\label{waveequationn}
\end{eqnarray}
where the stress tensor $\sigma$ is given by
\begin{eqnarray}
\sigma =\lambda I_d Tr(\mathcal{E}\,u^0)+2\mu \mathcal{E}u^0,
\label{stress}
\end{eqnarray}
$I_d$ is the identity on $\mathbb{R}^d$, and $Tr(\mathcal{E}\,u^0)$ is the trace of the strain (see \cite{Lipton2}). 
(See \cite{Lipton1} for a similar conclusion  associated with an alternative set of hypotheses.) 
The convergence of the peridynamic equation of motion to the local linear elastodynamic equation away from the crack set is 
consistent with the convergence of peridynamic equation of motion for {\emph{convex}} peridynamic potentials 
as seen in \cite{silling08}, \cite{Du-mengesha}, \cite{emmrich07}.

\begin{figure} 
\centering
\includegraphics[width=.6\textwidth]{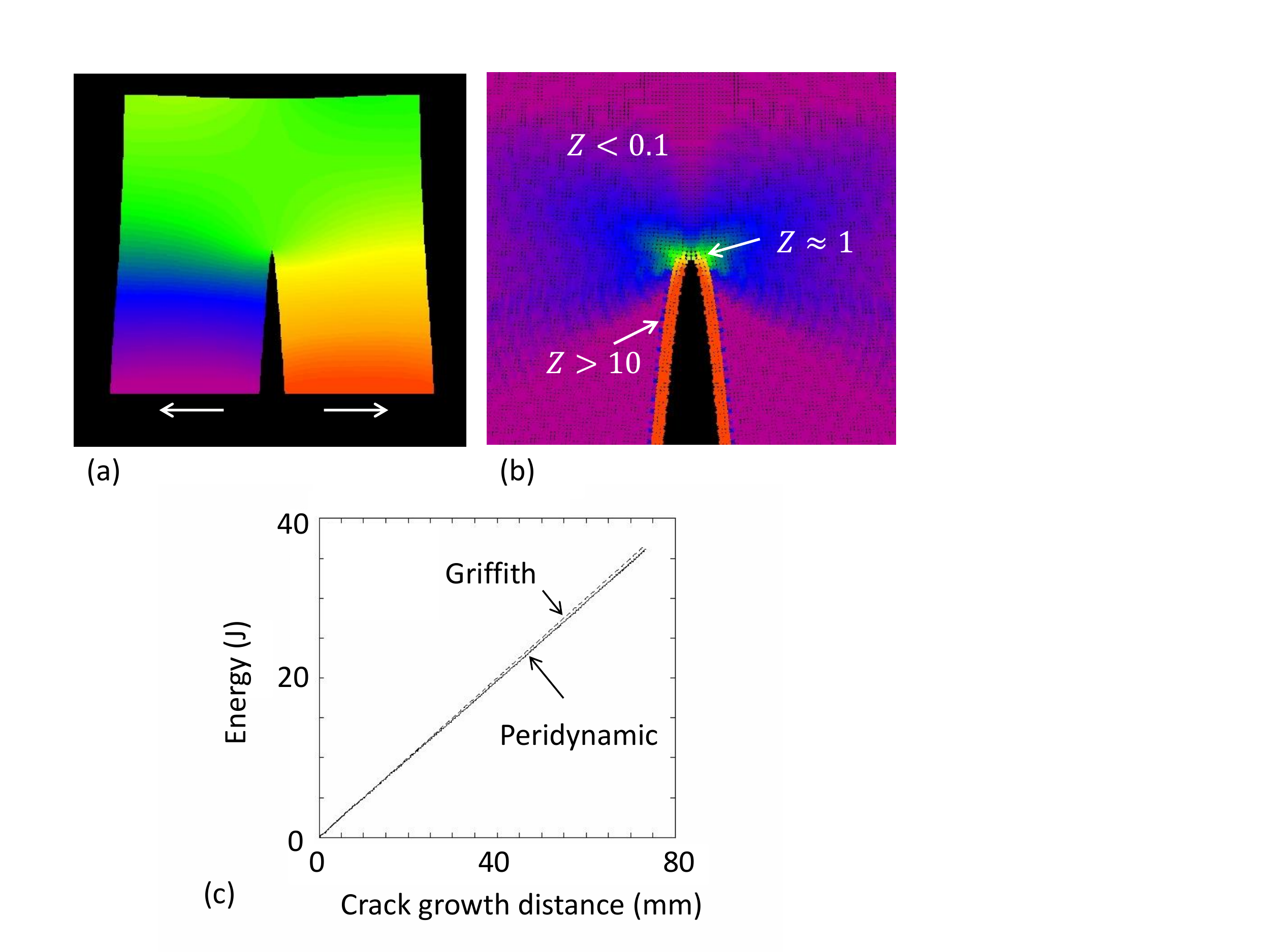}
\caption{Example 1: Stable dynamic crack growth.
(a) Displacement contours ($\times$100).
(b) Bond strain relative to the onset of instability near the crack tip.
(c) Energy consumed by the crack as a function of crack length.}
\label{fig-ex1}
\end{figure}

\section{Numerical Examples}

We present two example problems that demonstrate (1) that the nonlocal elastodynamic model reproduces a constant,
prescribed value of $\cG$, and (2) the model predicts reasonable behavior for the nucleation and propagation of complex
patterns of brittle dynamic fracture.

In the first example, a 0.1m $\times$ 0.1m plate with unit thickness has a material model of the form \eqref{eqn-WJPsi} 
with $\Psi(p) = c ( 1-e^{-\beta p})$ and $J(q)=1-q$
%\[   
%\Psi(p) = C \left( 1-e^{-\beta p}  \right), \quad  J(q)=1-\frac{q}{\eps}, 
%\]
where $c$ and $\beta$ are positive constants.
These constants are determined so that the bulk modulus $k$ and the critical energy release rate $\cG$ are
$k$=25GPa, $\cG$=500Jm$^{-2}$. 
The density is $\rho$=1200kg-m$^{-3}$
The maximum in the bond force curve occurs at $\cS_c=1/\sqrt{2\beta|y-x|}$.

An initial edge crack of length 0.02m extends vertically from the midpoint of the lower boundary (Figure~3).
A strip of thickness $\eps$ along the lower boundary is subjected to a constant velocity condition
$v_x=\pm 1.0$m/s, causing the crack to grow.
The solution method described in \cite{silling05} is used with a 400 $\times$ 400 square grid of nodes.
The horizon is $\eps$=0.00075m.
Figure~3(a) shows contours of displacement after the crack has grown  halfway through the plate.
The crack has a limiting growth velocity of about 1400m/s, which is about 50\% of the shear wave speed.
Figure~3(b) shows a close-up view of the growing crack tip  with
the colors indicating $Z(x)$. % where $Z(x)$ is the maximum value of bond strain relative
%to the critical strain $\cS_c$ among all bonds connected to $x$:
%\begin{equation}
%  Z(x)=\max_{y\in\cH_\eps(x)}\frac{\cS_u(x,y)}{\cS_c(x,y)}
%\label{eqn-Sdef}
%\end{equation}
Near the crack tip, the green lobes indicate a process zone in which the material goes through
a neutrally stable phase as bonds approach the maximum of the force vs. bond strain curve.

Figure~3(c) illustrates the energy balance in the model.
The curve labeled ``Griffith'' represents the idealized result under the assumption
that the crack uses a constant amount of energy $\cG$ per unit distance as it grows.
The curve labeled ``peridynamic'' is the energy that is stored in bonds in the computational model that have $Z\gg1$, that is, bonds
that are so far out on the curve in Figure~1 that they sustain negligible force density.
The fracture energy $\cG$ is stored in these bonds.
As shown in Figure~3(c), the energy consumed by the crack in the numerical model closely approximates
what is expected for a Griffith crack.
We conjecture that the small difference between the Griffith and peridynamic curves is due to numerical dissipation.

\begin{figure}  
\centering
\includegraphics[width=.45\textwidth]{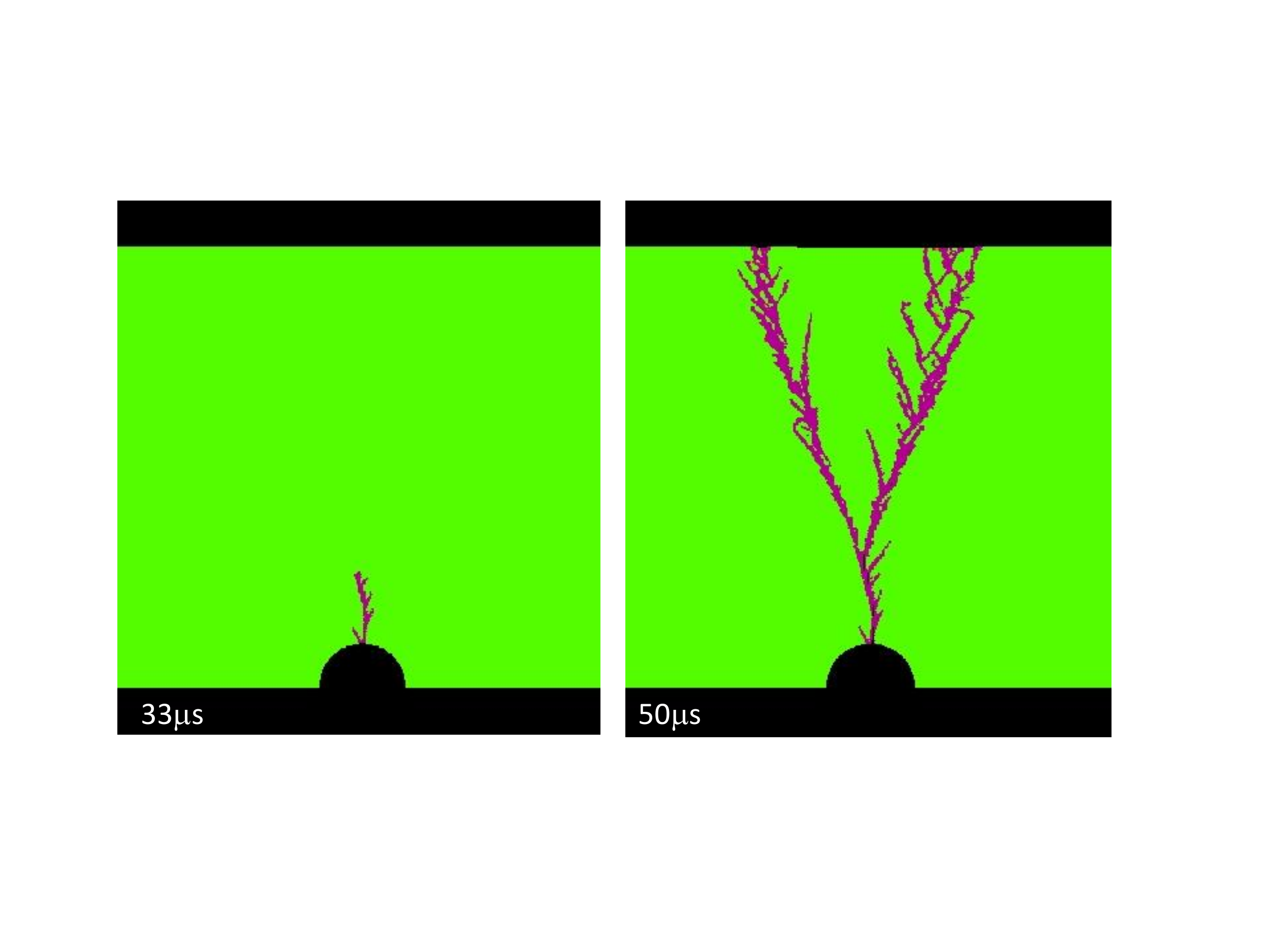}
\caption{Example 2: Computed paths of dynamic fractures nucleated at a circular notch
soon after nucleation (left) and after progression through the plate (right).}
\label{fig-ex2}
\end{figure}

In the second numerical example, the same material as above occupies a 0.2m $\times$ 0.1m rectangle in the
plane and has a semicircular notch as shown in Figure~4.
The material has an initial velocity field $v_1$=40m-s$^{-1}$, $v_2$=-13.3m-s$^{-1}$ throughout, where the $1$ and $2$
coordinates are in the horizonal and vertical directions, respectively.
The rectangular region has constant velocity boundary
conditions on the left and right boundaries that are consistent with the initial velocity field.
As time progresses, the strain concentration near the notch causes some bonds to exceed $Z=1$.
The resulting material instability nucleates cracks at the notch that rapidly accelerate and branch. The points $x$ associated $Z(x)>1$ are illustrated in  Figure~4  and correspond to the crack paths.
Many microbranches are visible in the crack paths.
For most of these microbranches, the strain energy is not sufficient to sustain growth, and they arrest.
Such microbranches are frequently seen in experiments on dynamic brittle fracture, for example \cite{fineberg1999}.

%\begin{figure} % [f]
%\includegraphics[width=.6\textwidth]{pnas-fig-ex1}
%\caption{Example 1: Stable dynamic crack growth.
%(a) Displacement contours ($\times$100).
%(b) Bond strain relative to the onset of instability near the crack tip.
%(c) Energy consumed by the crack as a function of crack length.}
%\label{fig-ex1}
%\end{figure}

%\begin{figure}  [f]
%\includegraphics[width=.45\textwidth]{pnas-fig-ex2}
%\caption{Example 2: Computed paths of dynamic fractures nucleated at a circular notch
%soon after nucleation (left) and after progression through the plate (right).}
%\label{fig-ex2}
%\end{figure}

%\newpage

\section{Observations and Discussion}

In this article we describe a theoretical and computational framework for analysis of complex brittle fracture based upon Newtons Second Law. 
This is enabled by recent advances in nonlocal continuum mechanics that treat singularities such as cracks according to the same
field equations and material model as points away from cracks.
This approach is different from  other contemporary approaches that involve the use of a phase field or 
cohesive zone elements to represent the fracture set, see \cite{bourdin,LarsenOrtnerSuli,Miehe,Borden,Duarte}. 

The key aspect of the elastic peridynamic material model that leads to crack growth is the nonconvexity 
of the bond energy density function.
In the classical theory of solid mechanics, nonconvex strain energy densities are related to the
emergence of features such as martensitic phase boundaries and crystal
twinning associated with the loss of ellipticity, a type of material instability.
As shown in the present paper, nonconvexity in peridynamic mechanics leads to crack nucleation and growth through an analogous material instability
within the nonlocal mathematical description.

%\appendix[Estimating the Spectral Norm of a Matrix]

%\appendix
%This is an example of an appendix without a title.

\section*{Acknowledgements}
This work was partially supported by 
NSF Grant DMS-1211066, AFOSR grant FA9550-05-0008, and 
NSF EPSCOR Cooperative Agreement No. EPS-1003897 with additional 
support from the Louisiana Board of Regents (to R.L.).
Sandia National Laboratories is a multi-program laboratory operated by Sandia Corporation, a wholly owned subsidiary of Lockheed Martin Corporation, 
for the U.S. Department of Energy's National Nuclear Security Administration under contract DE-AC04-94AL85000.

%\end{article}
\end{document}